\documentclass[12pt]{article}
\usepackage[all, knot]{xy}
\usepackage{xypic}
\usepackage{epsfig}
\usepackage{xspace}
\usepackage{layout}
\usepackage{latexsym}
\usepackage{amsmath}
\usepackage{amsthm}
\usepackage{amssymb}
\usepackage{amsfonts}
\usepackage{amsxtra}     
\usepackage{graphicx}
\usepackage{url}
\usepackage{verbatim}
\usepackage{longtable}

\input xypic
\input xy
\xyoption{all} \xyoption{arc}

\newcommand{\ie}{\text{i.e., }}
\newcommand{\eg}{\text{e.g., }}

\newcommand{\longdash}{\rule[.6ex]{3em}{.03ex}}
\newcommand{\Sym}{\text{\rm Sym}}

\begin{document}

\newtheorem{thm}{Theorem}[section]
\newtheorem{conj}{Conjecture}[section]
\newtheorem{lem}[thm]{Lemma}
\newtheorem{cor}[thm]{Corollary}
\newtheorem{prop}[thm]{Proposition}
\newtheorem{rem}[thm]{Remark}

\theoremstyle{definition}
\newtheorem{defn}[thm]{Definition}
\newtheorem{examp}[thm]{Example}
\newtheorem{notation}[thm]{Notation}
\newtheorem{rmk}[thm]{Remark}

\theoremstyle{remark}

\makeatletter
\renewcommand{\maketag@@@}[1]{\hbox{\m@th\normalsize\normalfont#1}}%
\makeatother
\renewcommand{\theenumi}{\roman{enumi}}

\def\square{\hfill${\vcenter{\vbox{\hrule height.4pt \hbox{\vrule width.4pt
height7pt \kern7pt \vrule width.4pt} \hrule height.4pt}}}$}

\newenvironment{pf}{{\it Proof:}\quad}{\square \vskip 12pt}

\title{Musical Actions of Dihedral Groups}
\author{Alissa S.~Crans, Thomas M.~Fiore, and Ramon Satyendra}


\maketitle

\section{Introduction}
Can you hear an action of a group? Or a centralizer? If knowledge of
group structures can influence how we {\it see} a crystal, perhaps
it can influence how we {\it hear} music as well. In this article we
explore how music may be interpreted in terms of the group structure
of the dihedral group of order 24 and its centralizer by explaining
two musical actions.\footnote{The composer Milton Babbitt was one of
the first to use group theory to analyze music. See
\cite{babbitt1960}.} The {\it dihedral group of order 24} is the
group of symmetries of a regular 12-gon, that is, of a 12-gon with
all sides of the same length and all angles of the same measure.
Algebraically, the dihedral group of order 24 is the group generated
by two elements, $s$ and $t$, subject to the three relations
$$s^{12}=1, \hspace{.75in} t^{2}=1, \hspace{.75in} tst=s^{-1}.$$

The first musical action of the dihedral group of order 24 we
consider arises via the familiar compositional techniques of {\it
transposition} and {\it inversion}. A transposition moves a sequence
of pitches up or down. When singers decide to sing a song in a
higher register, for example, they do this by transposing the
melody. An inversion, on the other hand, reflects a melody about a
fixed axis, just as the face of a clock can be reflected about the
0-6 axis. Often, musical inversion turns upward melodic motions into
downward melodic motions.\footnote {A precise, general definition of
inversion will be given later.} One can hear both transpositions and
inversions in many fugues, such as Bernstein's ``Cool'' fugue from
{\it West Side Story} or in Bach's {\it Art of Fugue}. We will
mathematically see that these musical transpositions and inversions
are the symmetries of the regular 12-gon.

The second action of the dihedral group of order 24 that we explore
has only come to the attention of music theorists in the past two
decades. Its origins lie in the $P,L,$ and $R$ operations of the
19th-century music theorist Hugo Riemann. We quickly define these
operations for musical readers now, and we will give a more detailed
mathematical definition in Section \ref{section:PLR}. The {\it
parallel} operation $P$ maps a major triad\footnote{A {\it triad} is
a three-note chord, \ie a set of three distinct pitch classes. {\it
Major} and {\it minor} triads, also called {\it consonant triads},
are characterized by their interval content and will be described in
Section \ref{section:majorminortriads}.} to its parallel minor and
vice versa. The {\it leading tone exchange} operation $L$ takes a
major triad to the minor triad obtained by lowering only the root
note by a semitone. The operation $L$ raises the fifth note of a
minor triad by a semitone. The {\it relative} operation $R$ maps a
major triad to its relative minor, and vice versa. For example,
$$P(C\text{-major})=c\text{-minor},$$
$$L(C\text{-major})=e\text{-minor},$$
$$R(C\text{-major})=a\text{-minor}.$$ It is through these three operations
$P,L,$ and $R$ that the dihedral group of order 24 acts on the set
of major and minor triads.

The $P,L,$ and $R$ operations have two beautiful geometric
presentations in terms of graphs that we will explain in Section
\ref{section:PLR}. Musical readers will quickly see that the
$C$-major triad shares two common tones with each of the three
consonant triads $P(C\text{-major})$, $L(C\text{-major})$, and
$R(C\text{-major})$ displayed above. These common tone relations are
geometrically presented by a toroidal graph with vertices the
consonant triads and with an edge between any two vertices having
two tones in common. This graph is pictured in two different ways in
Figures \ref{douthettsteinbach} and \ref{wallertorus}. As we shall
see, Beethoven's {\it Ninth Symphony} traces out a path on this
torus.\footnote{The interpretation of the {\it Ninth Symphony}
excerpt as a path on the torus was proposed by Cohn in
\cite{cohn1997}.}

Another geometric presentation of the $P,L$, and $R$ operations is
the {\it Tonnetz} graph pictured in Figure \ref{tonnetz}. It has
pitch classes as vertices and decomposes the torus into triangles.
The three vertices of any triangle form a consonant triad, and in
this way we can represent a consonant triad by a triangle. Whenever
two consonant triads share two common tones, the corresponding
triangles share the edge connecting those two tones. Since the
$P,L,$ and $R$ operations take a consonant triad to another one with
two notes in common, the $P,L,$ and $R$ operations correspond to
reflecting a triangle about one of its edges. The graph in Figures
\ref{douthettsteinbach} and \ref{wallertorus} is related to the {\it
Tonnetz} in Figure \ref{tonnetz}: they are {\it dual graphs}.

In summary, we have two ways in which the dihedral group acts on the
set of major and minor triads: (i) through applications of
transposition and inversion to the constituent pitch classes of any
triad, and (ii) through the operations $P,L$, and $R$. Most
interestingly, these two group actions are {\it dual} in the precise
sense of David Lewin \cite{lewin1987}. In this article we illustrate these group
actions and their duality in musical examples by Pachelbel, Wagner,
and Ives.

We will mathematically explain this duality in more detail later,
but we give a short description now.  First, we recall that the {\it
centralizer} of a subgroup $H$ in a group $G$ is the set of elements
of $G$ which commute with all elements of $H$, namely $$C_G(H)=\{g
\in G \mid gh=hg \text{ for all }h \in H \}.$$ The centralizer of
$H$ is itself a subgroup of $G$. We also recall that an action of a
group $K$ on a set $S$ can be equivalently described as a
homomorphism from $K$ into the symmetric group\footnote{The {\it
symmetric group} on a set $S$ consists of all bijections from $S$ to
$S$. The group operation is function composition.} $\Sym(S)$ on the
set $S$. Thus, each of our two group actions of the dihedral group
above gives rise to a homomorphism into the symmetric group on the
set $S$ of major and minor triads. It turns out that each of these
homomorphisms is an embedding, so that we have two distinguished
copies, $H_1$ and $H_2$, of the dihedral group of order 24 in
$\Sym(S)$. One of these copies is generated by $P,L$, and $R$. With
these notions in place, we can now express David Lewin's idea of
duality in \cite{lewin1987}: the two group actions are {\it dual} in
the sense that each of these subgroups $H_1$ and $H_2$ of $\Sym(S)$
is the centralizer of the other!

Practically no musical background is required to enjoy this
discussion since we provide mathematical descriptions of the
required musical notions, beginning with the traditional translation
of pitch classes into elements of $\mathbb{Z}_{12}$ via Figure
\ref{musicalclock}. From there we develop a musical model using
group actions and topology. We hope that this article will resonate
with mathematical and musical readers alike.

\section{Pitch Classes and Integers Modulo 12}

As the ancient Greeks noticed, any two pitches that differ by a
whole number of octaves\footnote{A pitch $y$ is an octave above a
pitch $x$ if the frequency of $y$ is twice that of $x$.} sound
alike. Thus we identify any two such pitches, and speak of {\it
pitch classes} arising from this equivalence relation. Like most
modern music theorists, we use equal tempered tuning, so that the
octave is divided into twelve pitch classes as follows.
$$\begin{tabular}{|c|c|c|c|c|c|c|c|c|c|c|c|c|}
\hline
$A$ & $A \sharp$ & $B$ & $C$ & $C \sharp$ & $D$ & $D \sharp$ & $E$ & $F$ & $F \sharp$ & $G$ & $G \sharp$ & $A$
\\  & $B \flat$ & & & $D \flat$  & & $E \flat$ &  &  & $G \flat$ & &
$A \flat$ & \\
\hline
\end{tabular}$$
The interval between two consecutive pitch classes is called a {\it
half-step} or {\it semitone}. The notation $\sharp$ means to move up
a semitone, while the notation $\flat$ means to move down a
semitone. Note that some pitches have two letter names. This is an
instance of {\it enharmonic equivalence}.

\begin{figure}[h]
\begin{center}
\includegraphics[height=2.5in]{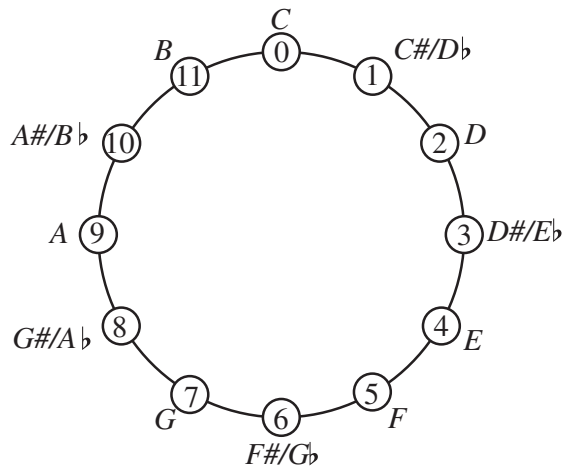}
\caption{The musical clock.}\label{musicalclock}
\end{center}
\end{figure}

Music theorists have found it useful to translate {\it pitch
classes} to {\it integers modulo 12} taking $0$ to be $C$ as in
Figure \ref{musicalclock}.
Mod 12 addition and subtraction can be read off of this clock; for
example $2+3=5$ mod 12, $11+4=3$ mod 12, and $1-4=9$ mod 12. We can
also determine the musical interval from one pitch class to another;
for example, the interval from $D$ to $G\sharp$ is six semitones.
This description of pitch classes in terms of $\mathbb{Z}_{12}$ can
be found in many articles, such as \cite{mccartin1998} and
\cite{rahn1980}. This translation from pitch classes to integers
modulo 12 permits us to easily use abstract algebra for modeling
musical events,  as we shall see in the next two sections.

\section{Transposition and Inversion}
Throughout the ages, composers have drawn on the musical tools of
transposition and inversion. For example, we may consider a type of
musical composition popular in the 18th century that is especially
associated with J.~S.~Bach: the {\it fugue}. Such a composition
contains a principal melody known as the {\it subject}; as the fugue
progresses, the subject typically will recur in transposed and
inverted forms. Mathematically speaking, {\it transposition} by an
integer $n$ mod 12 is the function
$$\xymatrix{T_n:\mathbb{Z}_{12} \ar[r] & \mathbb{Z}_{12}}$$
$$T_n(x):=x+n \text{ mod 12}$$
and {\it inversion}\footnote{At this point in our discussion,
musically experienced readers may notice that the word {\it
inversion} has several meanings in music theory.  The kind of
inversion we define here is different from {\it chord inversion} in
which pitches other than the root are placed in the bass.  This
latter kind of inversion accounts for terms such as {\it
first-inversion triad}.  Our discussion is not concerned with chord
inversion.} about $n$ is the function
$$\xymatrix{I_n:\mathbb{Z}_{12} \ar[r] & \mathbb{Z}_{12}}$$
$$I_n(x):=-x+n\text{ mod 12}.$$ Bach often used diatonic transposition and inversion, which we can
view as mod 7 transposition and inversion after identifying the
diatonic scale with $\mathbb{Z}_7$. However, many contemporary
composers intensively use mod 12 transposition and inversion; see
for example \cite{forte1977}, \cite{morris1988}, and
\cite{rahn1980}.

As is well known, these transpositions and inversions have a
particularly nice representation in terms of the musical clock in
Figure \ref{musicalclock}. The transposition $T_1$ corresponds to
clockwise {\it rotation} of the clock by $\frac{1}{12}$ of a turn,
while $I_0$ corresponds to a {\it reflection} of the clock about the
0-6 axis. Hence $T_1$ and $I_0$ generate the {\it dihedral group} of
symmetries of the 12-gon. Since $(T_1)^n=T_n$ and $T_n \circ
I_0=I_n$, we see that the 12 transpositions and 12 inversions form
the dihedral group of order 24. The compositions
$$T_m \circ T_n=T_{m+n\text{ mod 12}}$$
$$T_m \circ I_n=I_{m+n\text{ mod 12}}$$
$$I_m \circ T_n=I_{m-n\text{ mod 12}}$$
$$I_m \circ I_n=T_{m-n\text{ mod 12}}$$
are easy to verify. This group is often called the $T/I${\it
-group}. The first action of the dihedral group of order 24 on the
set of major and minor triads that we study is defined via the
$T/I$-group.

\section{Major and Minor Triads} \label{section:majorminortriads}

Triadic harmony has been in use for hundreds of years and is still
used every day in popular music.  In this section we use the
integers modulo 12 to define major and minor triads; in this way we
can consider them as objects upon which the dihedral group of order
24 may act.

\begin{figure}
\begin{center}
\includegraphics[height=2.5in]{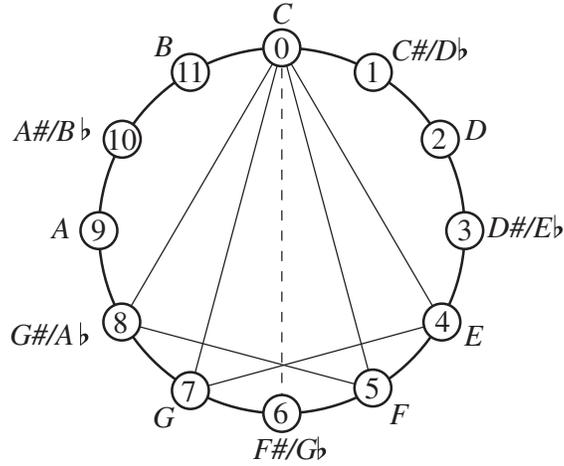}
\caption{$I_0$ applied to a $C$-major triad yields an $f$-minor
triad.}\label{cmajortriad}
\end{center}
\end{figure}

A {\it triad} consists of three simultaneously played notes.  A {\it
major triad} consists of a {\it root} note, a second note 4
semitones above the root, and a third note 7 semitones above the
root. For example, the $C$-major triad consists of
$\{0,4,7\}=\{C,E,G\}$ and is represented as a chord polygon in
Figure \ref{cmajortriad}. See \cite{mccartin1998} for beautiful
illustrations of the utility of chord polygons. Since any major
triad is a subset of the pitch-class space $\mathbb{Z}_{12}$, and
transpositions and inversions act on $\mathbb{Z}_{12}$, we can also
apply transpositions and inversions to any major triad. Figure
\ref{cmajortriad} shows what happens when we apply $I_0$ to the
$C$-major triad. The resulting triad is {\it not} a major triad, but
instead a {\it minor} triad.

A {\it minor triad} consists of a {\it root} note, a second note 3
semitones above the root, and a third note 7 semitones above the
root. For example, the $f$-minor triad consists of
$\{5,8,0\}=\{F,A\flat,C\}$ and its chord polygon appears in Figure
\ref{cmajortriad}.

Altogether, the major and minor triads form the set $S$ of {\it
consonant triads}, which are called {\it consonant} because of their
smooth sound.  A consonant triad is named after its root. For
example, the $C$-major triad consists of $\{0,4,7\}=\{C,E,G\}$ and
the $f$-minor triad consists of $\{5,8,0\}=\{F,A\flat,C\}$.
Musicians commonly denote major triads by upper-case letters and
minor triads by lower-case letters as indicated in the table of all
consonant triads in Figure \ref{majorminor}.
\begin{figure}
$$\begin{tabular}{|r|l|}
\hline \text{Major Triads} & \text{Minor Triads} \\ \hline
$C=\langle 0,4, 7\rangle$ & $\langle 0,8,5 \rangle = f $ \\
$C \sharp = D \flat=\langle 1,5, 8\rangle$ & $\langle 1,9, 6 \rangle=f \sharp=g \flat$ \\
$D=\langle 2,6, 9\rangle$ & $\langle 2,10,7\rangle = g$ \\
$D \sharp=E \flat=\langle 3,7, 10\rangle$ & $\langle3, 11,8 \rangle=g \sharp=a \flat$ \\
$E=\langle 4,8, 11\rangle$ & $\langle 4,0, 9\rangle=a$ \\
$F=\langle 5,9,  0 \rangle$ & $\langle 5,1,10  \rangle=a \sharp=b \flat$ \\
$F \sharp=G \flat=\langle 6,10, 1 \rangle$ & $\langle 6,2,11\rangle=b$ \\
$G=\langle 7,11, 2 \rangle$ & $\langle 7,3,0\rangle =c$ \\
$G \sharp=A \flat=\langle 8,0, 3 \rangle$ & $\langle8, 4,1\rangle =c \sharp=d \flat$ \\
$A=\langle 9,1, 4 \rangle$ & $\langle9, 5,2\rangle =d$ \\
$A \sharp=B \flat=\langle 10,2,5 \rangle$ & $\langle 10,6,3\rangle=d \sharp=e \flat$ \\
$B =\langle 11,3,6 \rangle$ & $\langle 11, 7,4\rangle =e$ \\
\hline
\end{tabular}$$
\caption{The set $S$ of consonant triads.}\label{majorminor}
\end{figure}

This table has several features. Angular brackets denote ordered
sets, which are called {\it pitch-class segments} in the music
literature. Since we are speaking of simultaneously sounding notes,
it is not necessary to insist on a particular ordering of the
elements within the brackets.\footnote{Another reason not to insist
on the ordering is the fact that the pitch-class set $\{0,4,7\}$ is
neither transpositionally nor inversionally symmetrical.} However
the mathematical artifice of an ordering will simplify the
discussion of the $PLR$-group and duality that we are approaching.
Such subtleties are discussed in \cite{fioresatyendra2005}.

The table also reflects the componentwise action of the $T/I$-group
because of this ordering. In the table, an application of $T_1$ to
an entry gives the entry immediately below it, for example
$$\aligned
T_1 \langle 0,4,7\rangle  &= \langle T_1(0),T_1(4),T_1(7)\rangle \\ &= \langle 1,5,8\rangle.
\endaligned$$ More generally, if we count the first entry as entry $0$, the $n$th entry in the first column is
\begin{equation} \label{transitivity1}
T_n\langle 0,4,7\rangle=\langle T_n(0),T_n(4),T_n(7) \rangle
\end{equation} and the $n$th entry in the second column is
\begin{equation} \label{transitivity2}
I_n\langle 0,4, 7\rangle=\langle I_n(0),I_n(4), I_n(7) \rangle.
\end{equation}

From the table we conclude that the action of the $T/I$-group is
{\it simply transitive}, that is, for any consonant triads $Y$ and
$Z$ there is a unique element $g$ of the $T/I$-group such that
$gY=Z$. As we have just seen in equations (\ref{transitivity1}) and
(\ref{transitivity2}), for any $Y$ and $Z$ there exist $g_1$ and
$g_2$ such that $g_1C=Z$ and $g_2C=Y$, and thus $gY=Z$ for
$g=g_1g_2^{-1}$. A quick verification also shows that $g$ is unique.

We can see the uniqueness of $g$ in a more elegant way using the
orbit-stabilizer theorem. The {\it orbit} of an element $Y$ of a set
$S$ under a group action of $G$ on $S$ consists of all those
elements of $S$ to which $Y$ is moved, in other words
$$\text{orbit of $Y$}=\{hY \mid h \in G \}.$$
The {\it stabilizer group} of $Y$ consists of all those elements of
$G$ which fix $Y$, namely
$$G_Y=\{h \in G \mid hY=Y\}.$$
\begin{thm}[Orbit-Stabilizer Theorem] \label{orbitstabilizer}
If a group $G$ acts on a set $S$ and $G_Y$ denotes the stabilizer group of $Y \in S$, then
$$|G|/|G_Y|=|\text{\rm orbit of $Y$}|.$$
\end{thm}
In our situation, $G$ is the dihedral group of order 24, $S$ is the
set of consonant triads as in Figure \ref{majorminor}, and
$|\text{orbit of $Y$}|=24$, so that $|G_Y|$=1. Thus, if $g'Y=gY$
then $g^{-1}g'Y=Y$, so that $g^{-1}g'$ is the identity element of
the group, and finally $g'=g$.

Generally, a group action of $G$ on a set $S$ is the same as a
homomorphism from $G$ into the symmetric group on the set $S$.
Indeed, from a group action we obtain such a homomorphism by
$$g \mapsto (Y \mapsto gY).$$ In the case of the $T/I$-group, this homomorphism is given by the
componentwise action of the $T/I$-group and it is injective.  For
simplicity we identify the $T/I$-group with its image in the
symmetric group on the set $S$.

\section{The $PLR$-Group} \label{section:PLR}
Up to this point, we have studied the action of the dihedral group
of order 24 on the set $S$ of major and minor triads via
transposition and inversion. Next we discuss a second musical action
of the dihedral group, but this time defined in terms of the
$PLR$-group.

Late 19th-century chromatic music, such as the music of Wagner, has
triadic elements to it but is not entirely tonal. For this reason,
it has been called ``triadic post-tonal'' in texts such as
\cite{cohn1998}. Recognizing that this repertoire has features which
are beyond the reach of traditional tonal theory, some music
theorists have worked on developing an alternative theory.

{\it Neo-Riemannian theory}, initiated by David Lewin in
\cite{lewin1982} and \cite{lewin1987}, has taken up the study of
$PLR$-transformations to address analytical problems raised by this
repertoire. We next define the {\it $PLR$-group} as the subgroup of
the symmetric group on the set $S$ generated by the bijections
$P,L,$ and $R$. As it turns out, this subgroup is isomorphic to the
dihedral group of order 24, as we prove in Theorem
\ref{PLRdihedral}. The $PLR$-group has a beautiful geometric
depiction in terms of a tiling on the torus called the {\it Tonnetz}
(Figure \ref{tonnetz}), which we also describe. A famous example
from Beethoven's {\it Ninth Symphony} is a path in the dual graph
(Figures \ref{douthettsteinbach} and \ref{wallertorus}).

\begin{figure}
\begin{center}
\noindent
\scalebox{.9}{\includegraphics[height=6.5in]{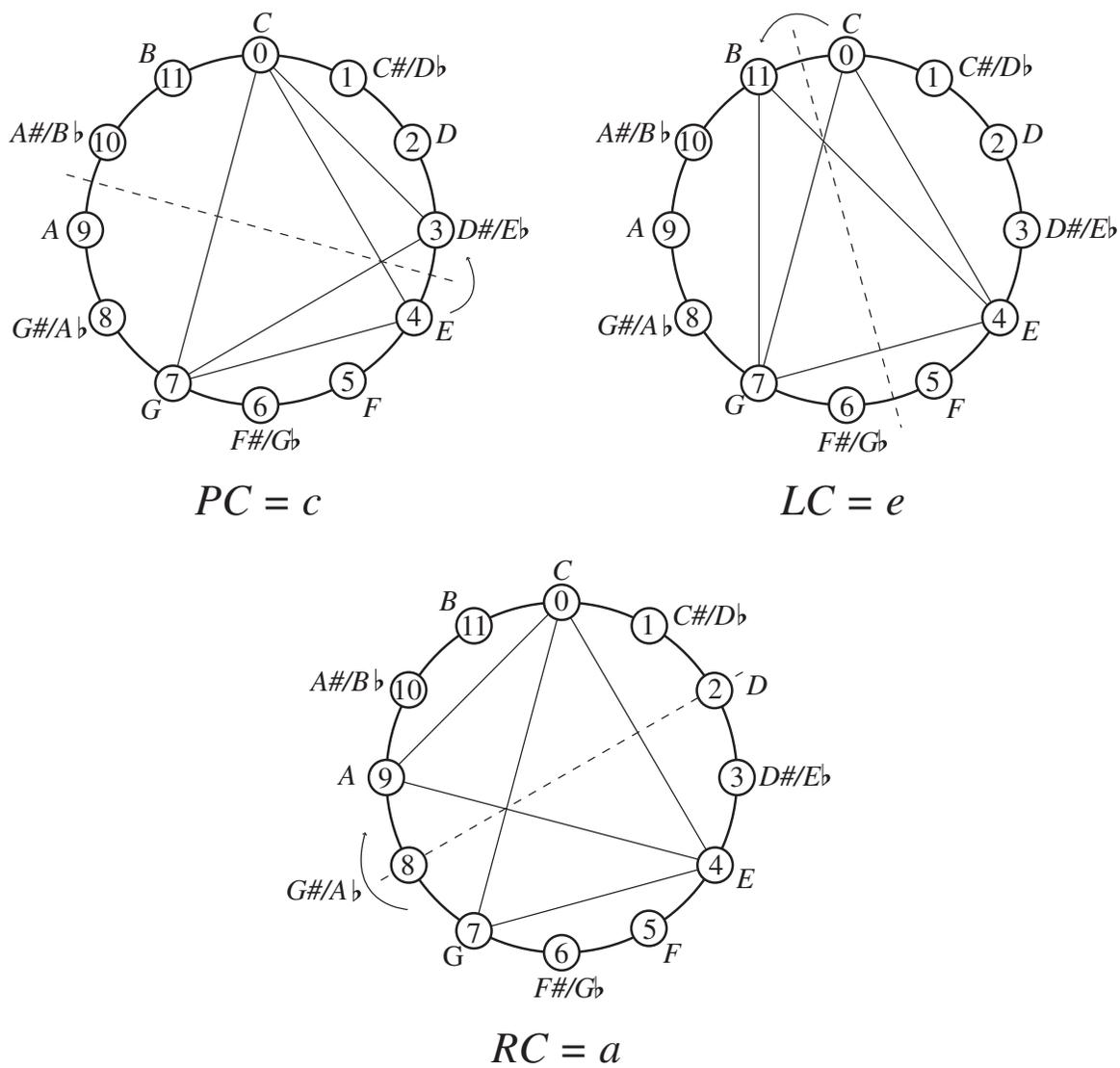}}
\caption{Minimal motion of the moving voice under $P,L,$ and
$R$.}\label{parsimony}
\end{center}
\end{figure}

\begin{figure}
\begin{center}
\includegraphics[height=3in]{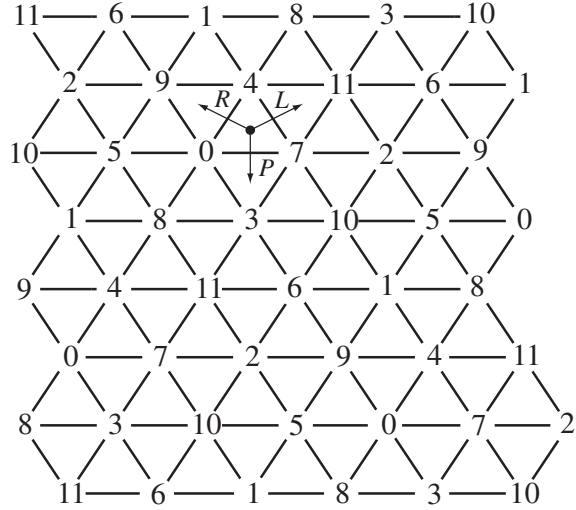}
\caption{The Oettingen/Riemann {\it Tonnetz}.} \label{tonnetz}
\end{center}
\end{figure}

\begin{figure}
\begin{center}
\includegraphics[height=2in]{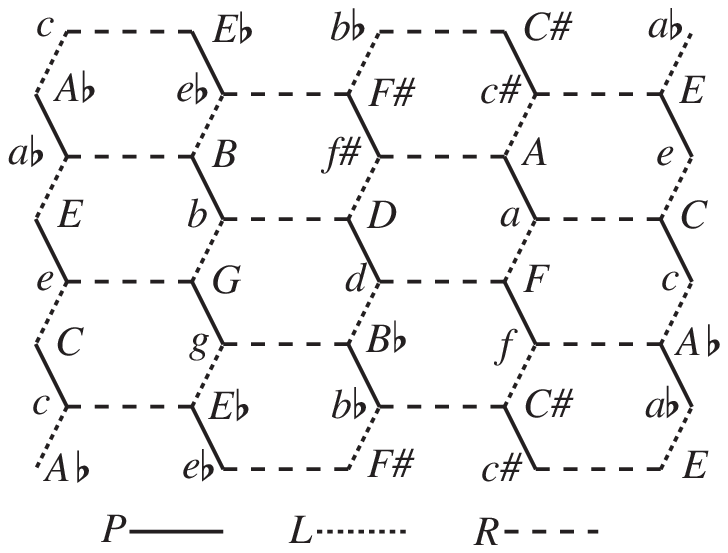}
\caption{Douthett and Steinbach's graph from
\cite{douthettsteinbach1998}.} \label{douthettsteinbach}
\end{center}
\end{figure}

\begin{figure}
\begin{center}
\includegraphics[height=6in]{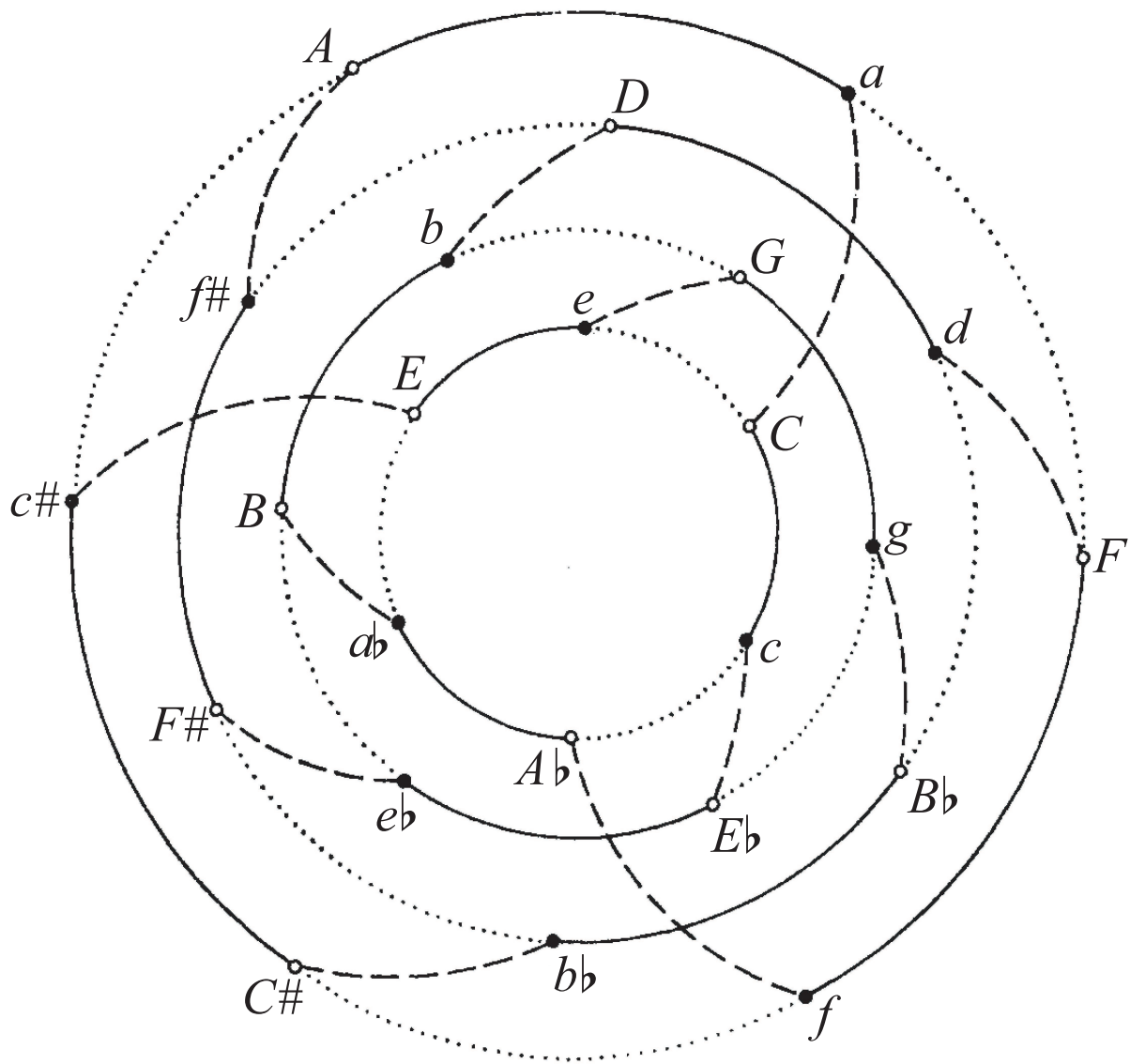}
\caption{Waller's torus from \cite{waller}.} \label{wallertorus}
\end{center}
\end{figure}

Consider the three functions $P,L,R:S \to S$ defined by
\begin{equation} \label{Palgebraic}
P\langle y_1, y_2, y_3\rangle=I_{y_1 + y_3}\langle y_1, y_2,
y_3\rangle
\end{equation}
\begin{equation} \label{Lalgebraic}
L\langle y_1, y_2, y_3\rangle=I_{y_2 + y_3}\langle y_1, y_2,
y_3\rangle
\end{equation}
\begin{equation} \label{Ralgebraic}
R\langle y_1, y_2, y_3\rangle=I_{y_1 + y_2}\langle y_1, y_2,
y_3\rangle.
\end{equation}
These are called {\it parallel, leading tone exchange,} and {\it
relative}. These are {\it contextual inversions} because the axis of
inversion depends on the aggregate input triad. Notably, the
functions $P,L,$ and $R$ are {\it not} defined componentwise, and
this distinguishes them from inversions of the form $I_n$, where the
axis of inversion is independent of the input triad. For $P,L,$ and
$R$ the axis of inversion on the musical clock when applied to
$\langle y_1, y_2, y_3\rangle$ is indicated in the table below.
\begin{center}
\begin{tabular}{|c|c|}
\hline
Function & Axis of Inversion Spanned by \\
$P$ & $\frac{y_1+y_3}{2},\frac{y_1+y_3}{2}+6$ \\
$L$ & $\frac{y_2+y_3}{2},\frac{y_2+y_3}{2}+6$ \\
$R$ & $\frac{y_1+y_2}{2},\frac{y_1+y_2}{2}+6$ \\
\hline
\end{tabular}
\end{center}
See Figure \ref{parsimony} for the axes of inversion in the
application of $P,L,$ and $R$ to the $C$-major triad.

If we consider major and minor as a {\it parity}, then there is a
particularly nice verbal description of $P,L,$ and $R$. The function
$P$ takes a consonant triad to that unique consonant triad of
opposite parity which has the first component and third component
switched. Thus, as unordered sets, the input and output triads
overlap in two notes. For example, $P\langle 0,4,7 \rangle=\langle
7,3,0 \rangle$ and $P\langle 7,3,0 \rangle=\langle 0,4,7 \rangle$. A
musician will notice that $P$ applied to $C$ is $c$, while $P$
applied to $c$ is $C$. In general, $P$ takes a major triad to its
parallel minor and a minor triad to its parallel major. A major
triad and a minor triad are said to be {\it parallel} if they have
the same letter name but are of opposite parity. The function $P$ is
manifestly an involution.

The other two functions, $L$ and $R$, similarly have maximally
overlapping inputs and outputs and are  involutions. The function
$L$ takes a consonant triad to that unique consonant triad of
opposite parity which has the second component and third component
switched; for example $L \langle 0,4,7 \rangle=\langle 11,7,4
\rangle$ and $L \langle 11,7,4 \rangle=\langle 0,4,7 \rangle $. The
function $R$ takes a consonant triad to that unique consonant triad
of opposite parity which has the first component and second
component switched; for example $R \langle 0,4,7 \rangle=\langle
4,0,9 \rangle$ and $R \langle 4,0,9 \rangle=\langle 0,4,7 \rangle$.
A musician will notice that $R$ applied to $C$ is $a$ and $R$
applied to $a$ is $C$. In general, $R$ takes a major triad to its
relative minor and a minor triad to its relative major. A major
triad and a minor triad are said to be {\it relative} if the root of
the minor triad is three semitones below the root of major triad.
The functions $R$ and $L$ are also involutions.

Each of the three functions corresponds to ubiquitous musical
motions that are easy to learn to recognize by ear.  That the input
and output overlap in two common tones is one reason the motions are
easily recognized. These three triadic transformations were employed
by European composers with great success in the years 1500-1900.
Another distinguishing feature is the minimal motion of the moving
voice. For example, in the application of these three functions to
the $C$-major triad above, we see in the case of $P$ that 4 moves to
3, in the case of $L$ that 0 moves to 11, and in the case of $R$
that 7 moves to 9. This is illustrated in Figure \ref{parsimony}.

This {\it parsimonious voice leading} is unique to the major and
minor triads as shown in \cite{cohn1997}: if one starts with any
other three note chord, such as $\langle 0,1,3 \rangle$ for example,
and generates 24 chords by transposition and inversion, then the
analogues of $P,L,$ and $R$ will always have large jumps in their
moving voices.\footnote{If one starts with $\langle 0,4,8 \rangle$,
then $P,L,$ and $R$ will be trivial, so we exclude this case.} As
Cohn points out in \cite{cohn1997}, the potential for parsimonious
voice leading is totally independent of the acoustic properties of
consonant triads; instead it is ``a function of their
group-theoretic properties as equally tempered entities modulo 12.''

The group generated by $P,L,$ and $R$ is called the {\it
$PLR$-group} or the {\it neo-Riemannian group} after the late
19th-century music theorist Hugo Riemann. Its structure is well
known, as we illustrate in the following theorem. An important
ingredient for our proof is a famous chord progression in
Beethoven's {\it Ninth Symphony}. Cohn observed this chord
progression in \cite{cohn1997}.

\begin{thm} \label{PLRdihedral}
The $PLR$-group is generated by $L$ and $R$ and is dihedral of order 24.
\end{thm}
\begin{pf}
First we remark that one can use formulas (\ref{Palgebraic}),
(\ref{Lalgebraic}), and (\ref{Ralgebraic}) to show that $PT_1=T_1P$,
$LT_1=T_1L$, and $RT_1=T_1R$.

If we begin with the $C$-major triad and alternately apply $R$ and
$L$, then we obtain the following sequence of triads.\footnote{We recall that upper-case letters
refer to major triads and lower-case letters refer to minor triads.}
$$C,a,F,d,B \flat, g, E
\flat, c, A \flat, f, D \flat, b \flat, G \flat, e \flat, B, g
\sharp, E, c \sharp, A, f \sharp, D, b, G, e, C$$ This tells us that
the 24 bijections $R,LR,RLR,\dots, R(LR)^{11},$ and $(LR)^{12}=1$
are distinct, that the $PLR$-group has at least 24 elements, and
that $LR$ has order 12. Further $R(LR)^3(C)=c$, and since $R(LR)^3$
has order 2 and commutes with $T_1$, we see that $R(LR)^3=P$, and
the $PLR$-group is generated by $L$ and $R$ alone.

If we set $s=LR$ and $t=L$, then $s^{12}=1, t^2=1,$ and
$$\aligned tst &= L(LR)L \\
&= RL \\
&= s^{-1}. \endaligned$$ It only remains to show that the
$PLR$-group has order 24, and then it will be dihedral as on page 68
of \cite{rotman1995}. We postpone the proof of this last fact until
Theorem \ref{maintheorem}.
\end{pf}
\begin{cor} \label{PLRdihedralcorollary}
The $PLR$-group acts simply transitively on the set of consonant
triads.
\end{cor}
\begin{pf}
From the chord progression in Theorem \ref{PLRdihedral} we see that
the orbit of $C$-major is all of $S$, and has 24 elements. As the
$PLR$-group also has 24 elements, simple transitivity follows from
the orbit-stabilizer theorem.
\end{pf}

The Oettingen/Riemann {\it Tonnetz} in Figure \ref{tonnetz} is a
beautiful geometric depiction of the $PLR$-group.
The word {\it Tonnetz} is German for ``tone network'' and is
sometimes translated as the ``table of tonal relations.'' The
vertices of this graph are pitch classes, while each of the
triangles is a major or minor triad. The graph extends infinitely in
all directions, though we have only drawn a finite portion. On the
horizontal axis we have the circle of fifths, and on the diagonal
axes we have the circles of major and minor thirds.\footnote{The
intervallic torus for minor thirds described in Table 2 of
\cite{mccartin1998} is contained in a diagonal of the {\it
Tonnetz}.} Since these circles repeat, we see that the {\it Tonnetz}
is doubly periodic. Therefore we obtain a torus by gluing the top
and bottom edges as well as the left and right edges of the
rectangular region indicated in Figure \ref{tonnetz}. The functions
$P,L,$ and $R$ allow us to navigate the {\it Tonnetz} by flipping a
triangle about an edge whose vertices are the preserved pitch
classes. This is investigated in \cite{cohn1997} for scales of
arbitrary chromatic number.

The Oettingen/Riemann {\it Tonnetz} in Figure \ref{tonnetz} is
similar to the one in Figure 2 on page 172 of
\cite{cohn1998}.\footnote{Our Figure 5 does not exactly reproduce
Figure 2 of \cite{cohn1998}, but introduces the following changes:
pitch-class numbers are shown rather than letter note names, the $D$
arrow is deleted, and a different region of the {\it Tonnetz} is
displayed. Special thanks go to Richard Cohn for giving us
permission to use this modified version of the figure.} Figure
\ref{tonnetz} is an {\it interpretation} of Riemann's {\it Tonnetz},
which resulted from the work of many neo-Riemannian theorists,
especially \cite{cohn1992}, \cite{hyer}, and
\cite{lewin1982}.\footnote{The article \cite{hyer} contains the
first appearance of the group generated by $P,L,R,$ and $D$, where
$D=T_5$ is the {\it dominant} transformation. This group appears
again in \cite{hooktriadic} as the group $\mathcal{H}$ on page 98.
Interestingly, $D=LR$ on major triads, but $D=RL$ on minor triads.}
Enharmonic equivalence and equal-tempered tuning are crucial for
this modern interpretation. Since Riemann did not use enharmonic
equivalence nor equal tempered tuning, his original {\it Tonnetz}
was not periodic and did not lie on a torus. The original {\it
Tonnetz} can be found on page 20 of \cite{riemannoriginal}, or on
page 102 of the translation \cite{riemanntranslation} (annotated in
\cite{wason}).

Douthett and Steinbach have associated the graph in Figure \ref{douthettsteinbach}
 to the neo-Riemannian $PLR$-group in \cite{douthettsteinbach1998}.\footnote{Figure
\ref{douthettsteinbach} has been reproduced by kind permission of
the authors.} This time the vertices are the consonant triads, and
there is an edge between two vertices labelled by $P,L,$ or $R$
whenever $P,L$, or $R$ sends one vertex to the other.  This graph is
also periodic vertically and horizontally, so the top and bottom
edges can be glued together, and the left and right edges can also
be glued after twisting a third of the way. The result is a graph on
the torus. Earlier, Waller studied this graph on the torus in
\cite{waller}, and observed that its automorphism group is the
dihedral group of order 24. Waller's torus is pictured in Figure
\ref{wallertorus}.\footnote{Waller's torus from \cite{waller} has
been reproduced in Figure 7 by kind permission of the U.K.
Mathematical Association and the {\it Mathematical Gazette}.}
Douthett and Steinbach also make this observation in
\cite{douthettsteinbach1998}, and present Waller's torus in the
context of neo-Riemannian theory.

Movement in music can be likened to movement along the surface of
the torus.  The sequence of consonant triads in the proof of Theorem
\ref{PLRdihedral} traces out a regular path on the torus in Figure
\ref{wallertorus}, and the first 19 triads of that sequence occur in
order in measures 143-176 of the second movement of Beethoven's {\it
Ninth Symphony}! Cohn observed this remarkable sequence in
\cite{cohn1992}, \cite{cohn1997}, and \cite{cohn1991}.

There is a relationship between the two graphs and their tori: they
are {\it dual graphs}. That means if you draw a vertex in the center
of every hexagonal face of Figure \ref{douthettsteinbach} or
\ref{wallertorus}, and connect two vertices by an edge whenever the
corresponding faces have a common edge, then you get the {\it
Tonnetz}.  In fact, a vertex of the {\it Tonnetz} is the unique note
in the intersection of the triads on the corresponding face; \eg 0
is the intersection of $a,C,c,A\flat,f,$ and $F$.

But in the musical model we are considering, these graphs are not
the only things which are dual. Using the notion of centralizer, we
will show that the $T/I$-group and the $PLR$-group are dual groups!

\section{$T/I$ and $PLR$ are Dual}
As we have seen, the dihedral group of order 24 acts on the set $S$
of major and minor triads simply transitively in two interesting
ways: (i) through the $T/I$-group using transposition and inversion,
and (ii) through the neo-Riemannian $PLR$-group using the $P,L,$ and
$R$ functions. If we consider the $T/I$-group and the $PLR$-group as
subgroups of the symmetric group $\Sym(S)$ on the set $S$, then an
interesting relation occurs: the centralizer of the $T/I$-group is
the $PLR$-group and the centralizer of the $PLR$-group is the
$T/I$-group! This means the $T/I$-group and the $PLR$-group are {\it
dual groups} in the terminology of Lewin \cite{lewin1987}. We prove
this momentarily. This duality in the sense of Lewin has also been
studied on pages 110-111 of \cite{hooktriadic}, and also in
\cite{hookthesis}.\footnote{In \cite{hooktriadic} and
\cite{hookthesis}, Hook embedded the neo-Riemannian $PLR$-group into
the group $\mathcal{U}$ of uniform triadic transformations. In the
following explanation of this embedding into Hook's group, we use
$S$ to denote the set of consonant triads, as in most of the present
article. A {\it uniform triadic transformation} $U$ is a function
$U:S \to S$ of the form $\langle\sigma,t^+,t^-\rangle$ where $\sigma
\in \{+,-\}$, and $t^+,t^- \in \mathbb{Z}_{12}$. The sign $\sigma$
indicates whether $U$ preserves or reverses parity (major
vs.~minor), the component $t^+$ indicates by how many semitones $U$
transposes the root of a major triad, and the component $t^-$
indicates by how many semitones $U$ transposes the root of a minor
triad. For example, the neo-Riemannian operation $R$ is written as
$\langle-,9,3\rangle$, meaning that $R$ maps any major triad to a
minor triad whose root is 9 semitones higher, and $R$ maps any minor
triad to a major triad whose root is 3 semitones higher, as one sees
with $R(C)=a$ and $R(a)=C$. Other familiar elements in $\mathcal{U}$
are $P=\langle-,0,0\rangle,$
 $L=\langle-,4,8\rangle,$ $R=\langle-,9,3\rangle,$ and $T_n=\langle+,n,n\rangle$. Uniform triadic
transformations are automatically invertible, like all these
examples. The non-Riemannian operations $D=T_5$ and
$M=\langle-,9,8\rangle$, called {\it dominant} and {\it diatonic
mediant} respectively, are also contained in $\mathcal{U}$. Thus,
the group $\mathcal{U}$ of uniform triadic transformations is a good
place to study how Riemannian operations and non-Riemannian
operations interact. However, the inversions $I_n$ are {\it not} in
$\mathcal{U}$. The uniform triadic transformations {\it and}
inversions are contained in the group $\mathcal{Q}$ of {\it quasi
uniform triadic transformations}. This group is much larger:
$|\mathcal{Q}|=1152$ while $|\mathcal{U}|=288$.

Hook defined on page 110 of \cite{hooktriadic} a {\it duality
operator} on $\mathcal{Q}$ which restricts to an anti-isomorphism
between the $T/I$-group and the $PLR$-group; transpositions and
inversions are mapped to {\it Schritte} and {\it Wechsel}
respectively. Morever, the Lewinnian duality we study in this paper
between $T/I$ and $PLR$ in $\Sym(S)$ restricts to the subgroup
$\mathcal{Q}$ of $\Sym(S)$: the centralizer of the $T/I$-group in
$\mathcal{Q}$ is precisely the $PLR$-group and the centralizer of
the $PLR$-group in $\mathcal{Q}$ is precisely the $T/I$-group.
Interestingly, the centralizer of the transposition group in
$\mathcal{Q}$ is $\mathcal{U}$. Even better, the centralizer of the
transposition group in $\Sym(S)$ is exactly $\mathcal{U}$ by Theorem
1.7 of \cite{hooktriadic}. The group $\mathcal{U}$ is isomorphic to
the wreath product $\mathbb{Z}_{12}\wr\mathbb{Z}_2$.}

The term ``dualism'' in the neo-Riemannian literature, such as
\cite{hooktriadic} and \cite{hookthesis},  is used mostly to refer
to a different idea associated with the music theorist Hugo Riemann.
Specifically, Riemannian ``dualism'' refers to a property of {\it
individual} elements of the $PLR$-group. A precise definition can be
found on page 59 of \cite{hooktriadic}: ``This property---whatever a
transformation does to a major triad, its effect on a minor triad is
precisely the opposite---may be regarded as an explicit
representation of Riemann's harmonic dualism.''

\begin{figure}
\begin{center}
\noindent\includegraphics[height=5.75in]{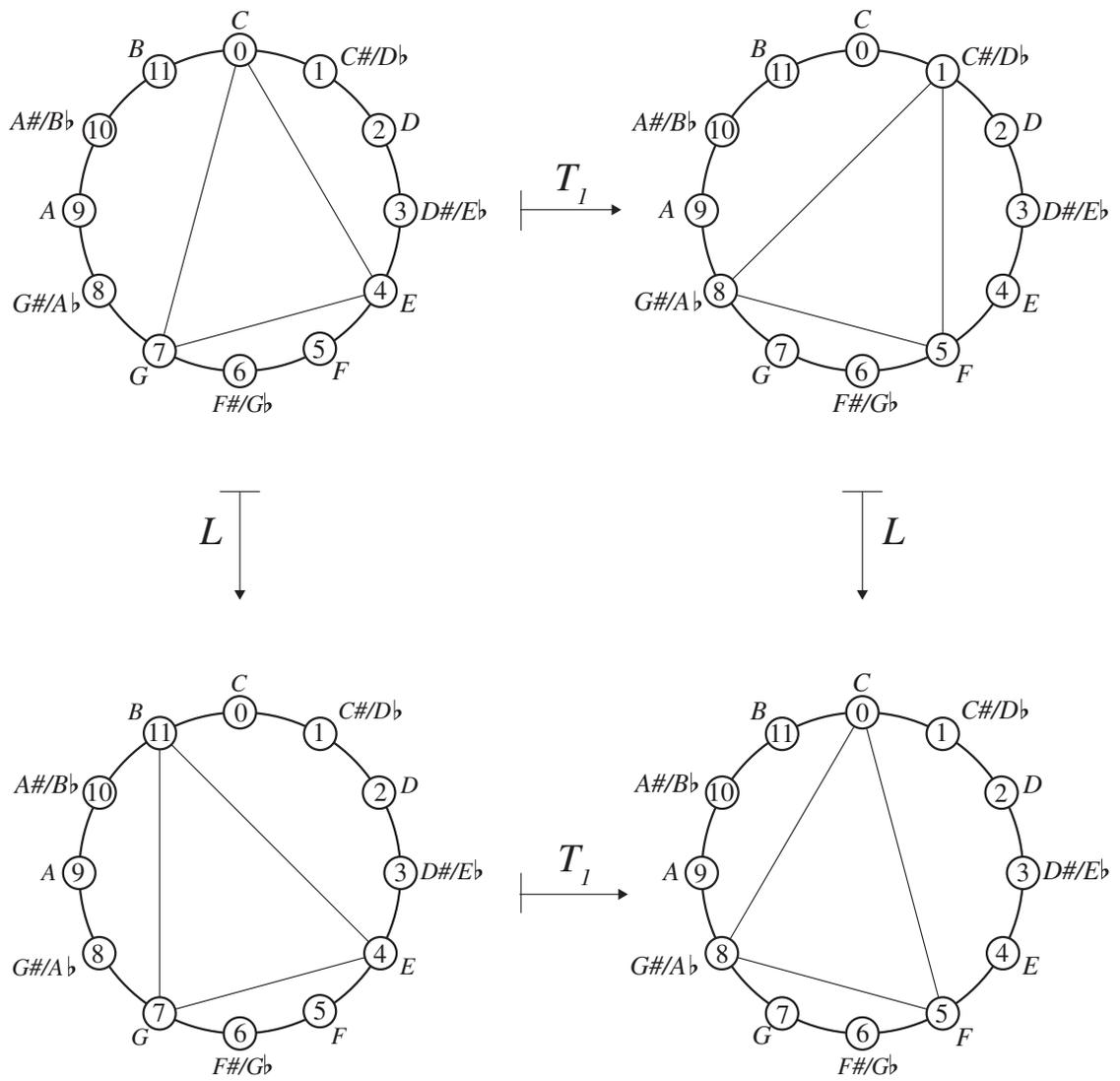}
\caption{Illustration of commutativity of $T_1$ and $L$.}
\label{commutativity}
\end{center}
\end{figure}

As an illustration of the duality between the $T/I$-group and the
$PLR$-group in the sense of Lewin, we can compute with the $C$-major
triad. If we apply $T_1$ to the $C$-major triad, and then $L$, that
is the same as first applying $L$ and then applying $T_1$ (see
Figure \ref{commutativity}). A category theorist would say that the
diagram
$$\xymatrix@R=3pc@C=3pc{S \ar[r]^{T_1} \ar[d]_L &  S \ar[d]^L \\
S \ar[r]_{T_1} & S}$$ {\it commutes}, \ie the result is the same no
matter which path one takes. Similarly, one can use formulas
(\ref{Palgebraic}), (\ref{Lalgebraic}), and (\ref{Ralgebraic}) to
show that $P,L,$ and $R$ commute with $T_1$ and $I_0$. Since these
are the generators of the respective groups, we conclude that any
diagram with vertical arrows in the $PLR$-group and horizontal
arrows in the $T/I$-group, as in Figure \ref{commutativity}, will
commute.

\begin{thm} \label{maintheorem}
The $PLR$-group and the $T/I$-group are dual. That is, each acts
simply transitively on the set $S$ of major and minor triads, and
each is the centralizer of the other in the symmetric group
$\Sym(S)$.
\end{thm}
\begin{pf}
In Section \ref{section:majorminortriads} we already concluded that
the $T/I$-group acts simply transitively on the set of major and
minor triads from Figure \ref{majorminor} and equations
(\ref{transitivity1}) and (\ref{transitivity2}). We also determined
in the discussion just before the statement of the current theorem
that any element of the $PLR$-group commutes with any element of the
$T/I$-group. In other words, the $PLR$-group is contained in the
centralizer $C(T/I)$ of the $T/I$-group in $\Sym(S)$.

For any element $Y$ of $S$ we claim that the stabilizer of $Y$ under
the action of $C(T/I)$ contains only the identity element. Suppose
that $h$ is in $C(T/I)$ and fixes $Y$, and that $g$ is in the
$T/I$-group. Then we have
$$\aligned
hY &= Y \\
ghY &= gY \\
hgY &= gY. \\
\endaligned$$
Since the $T/I$-group acts simply transitively, every $Y'$ in $S$ is
of the form $gY$ for some $g$ in the $T/I$-group, and therefore $h$
is the identity function on $S$ by the last equation above. Thus the
stabilizer $C(T/I)_Y$ of $Y$ in $C(T/I)$ is the trivial group.

An application of the orbit-stabilizer theorem to $G=C(T/I)$ gives
us
$$|C(T/I)|/|C(T/I)_Y|=|\text{\rm orbit of $Y$}|\leq |S|=24.$$
As the $PLR$-group is a subgroup of $C(T/I)$ and $|C(T/I)_Y|=1$, we
conclude $$|PLR\text{-group}| \leq |C(T/I)| \leq 24.$$

From the famous chord progression of Beethoven's {\it Ninth
Symphony} in the first part of Theorem \ref{PLRdihedral}, we know
that the $PLR$-group has at least 24 elements. Thus, the $PLR$-group
has exactly 24 elements and is equal to $C(T/I)$. This completes the
proof of Theorem \ref{PLRdihedral}, so we may now conclude as in
Corollary \ref{PLRdihedralcorollary} that the $PLR$-group acts
simply transitively on $S$.

It only remains to show that the $T/I$-group is the centralizer of
the $PLR$-group. However, this follows by reversing the roles of the
$T/I$-group and the $PLR$-group in the orbit-stabilizer argument we
just made.
\end{pf}

Now that we have met an important example of dual groups, one may
ask if there are other examples as well and how they arise.  Dual
groups have indeed been known for over 100 years, and can arise in
only one way, as the following theorem specifies.
\begin{thm}[Cayley]
If $G$ is a group, then we obtain dual groups via the two embeddings
of $G$ into $\Sym(G)$ as left and right actions of $G$ on itself.
All dual groups arise in this way.\footnote{We thank L\'{a}szl\'{o}
Babai for reminding us of this classical theorem.}
\end{thm}

\begin{figure}
\begin{center}
\includegraphics[width=3.5in]{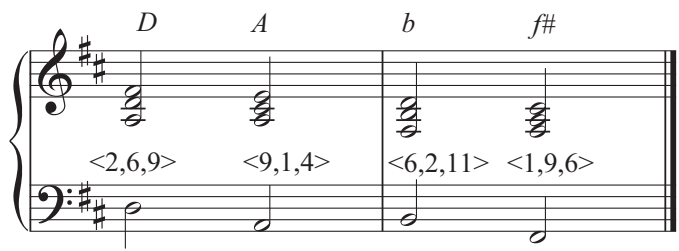}
\caption{Chord progression from Pachelbel, Canon in $D$.}
\label{Pachelbel}
\end{center}
\end{figure}

We now present three musical examples of the duality between the
$T/I$-group and the $PLR$-group. Our first example is Johann
Pachelbel's famous Canon in $D$, composed circa 1680 and reproduced
in Figure \ref{Pachelbel}.
The chord progression in the associated commutative diagram occurs
in 28 variations in the piece.
$$\xymatrix@R=4pc@C=4pc{D \ar@{|->}[r]^{T_7} \ar@{|->}[d]_R & A \ar@{|->}[d]^R \\ b \ar@{|->}[r]_{T_7} &
f\sharp}$$

\begin{figure}[h]
\begin{center}
\noindent \includegraphics[width=5in]{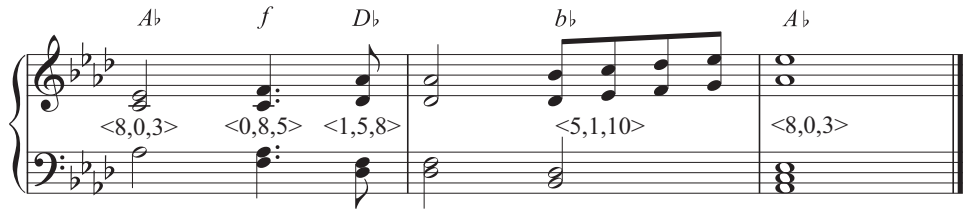}
\caption{Wagner, {\it Parsifal}, ``Grail'' Theme.} \label{Wagner}
\end{center}
\end{figure}

Another example can be found in the ``Grail'' theme of the Prelude
to Parsifal, Act 1, an opera completed by Richard Wagner in 1882.
See Figure \ref{Wagner} and the following commutative diagram.
$$\xymatrix@R=4pc@C=4pc{A\flat \ar@{|->}[d]_R \ar@{|->}[r]^{T_5} &  D\flat \ar@{|->}[d]^R
\\ f \ar@{|->}[r]_{T_5} & b \flat}$$

A particularly interesting example is in the opening measure of
``Religion,'' a song for voice and piano written by Charles Ives in
the 1920s. This time the horizontal transformation is an inversion,
namely $I_6$. Since the inversion $I_6$ transforms major triads to
minor triads, we have $LR$ acting upon triads of different parity.
This allows us to observe that $LR$ transforms $D$-major {\it up} by
5 semitones, but at the same time transforms $a$-minor {\it down} by
5 semitones. This makes the behavior of the left column {\it dual}
(in the sense of Riemann) to the behavior of the right column.
\begin{figure}[h]
\begin{center}
\includegraphics[width=3.5in]{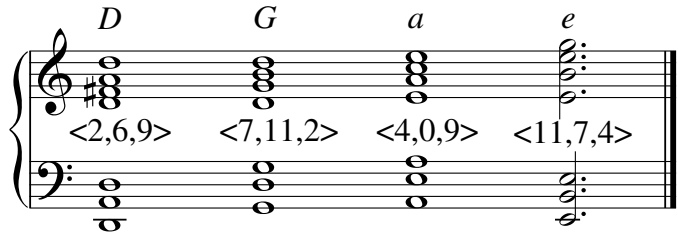}
\caption{Ives, ``Religion''.} \label{IVES}
\end{center}
\end{figure}
$$\xymatrix@R=4pc@C=4pc{D \ar@{|->}[d]_{LR} \ar@{|->}[r]^{I_6} &  a  \ar@{|->}[d]^{LR}
\\ G \ar@{|->}[r]_{I_6} & e}$$

\section{Recapitulation and Variation}
In summary, the dihedral group of order 24 acts on the set of major
and minor triads in two ways: through the $T/I$-group and through
the $PLR$-group. Further, these two actions are dual. The
$PLR$-group has two interesting geometric depictions: the {\it
Tonnetz} and Waller's torus. But why stop at major and minor triads?
One could just as well study the analogues of $P,L,$ and $R$ in the
context of dominant seventh chords and half-diminished seventh
chords. Indeed, that has been pursued in \cite{childs1998} and
\cite{gollin1998}. Moreover, the theory can be generalized further;
the authors of \cite{fioresatyendra2005} studied a neo-Riemannian
group for arbitrary pitch-class segments in terms of contextual
inversion, and applied their findings to an analysis of Hindemith,
{\it Ludus Tonalis}, Fugue in $E$. Neo-Riemannian groups for
asymmetrical pitch-class segments were studied in \cite{hooktriadic}
and \cite{hookthesis} from a root-interval point of view.

There are many avenues of exploration for undergraduates. Students
can listen to group actions in action and apply the orbit-stabilizer
theorem to works of music.

By experimenting with the $PLR$-group, students can also learn about
generators and relations for groups. The torus for Beethoven's {\it
Ninth Symphony} is an inviting way to introduce students to
topology. More tips for undergraduate study can be found on the
website \cite{fiorewebsite}, which contains lecture notes, problems
for students, slides, and more examples. For both advanced readers
and students, the website \cite{baezwebsite} includes entertaining
discussion and interesting posts by musicians and mathematicians
alike.

\paragraph{Acknowledgments.} Thomas M.~Fiore was supported at the University of Chicago
by NSF Grant DMS-0501208. At the Universitat Aut\`{o}noma de
Barcelona he was supported by Grant SB2006-0085 of the Programa
Nacional de ayudas para la movilidad de profesores de universidad e
investigadores espa$\tilde{\text{n}}$oles y extranjeros. Thomas
M.~Fiore and Alissa S.~Crans both thank Peter May for his ongoing
support and encouragement. Thomas M.~Fiore thanks Karen Rhea at the
University of Michigan for offering him the opportunity to try this
material in the classroom. He also thanks the students of Math 107
at the University of Michigan and the VIGRE REU students at the
University of Chicago who eagerly explored these ideas in class and on their own. The authors also thank Blake Mellor for assistance on the figures, and
Michael Orrison for comments on an earlier draft.


\bigskip

\noindent\textbf{Alissa S.~Crans} earned her B.S. in Mathematics
from the University of Redlands and her Ph.D. in Mathematics from
the University of California at Riverside, under the guidance of
John Baez.  She is currently an Assistant Professor of Mathematics
at Loyola Marymount University and has taught at Pomona College, the
Ohio State University, and the University of Chicago. Along with
Naiomi Cameron and Kendra Killpatrick, Alissa organizes the Pacific
Coast Undergraduate Mathematics Conference. In addition to
mathematics, she enjoys playing the clarinet with the LA Winds,
running, biking, reading, and traveling.

\noindent\textit{Department of Mathematics \\ Loyola Marymount
University \\ One LMU Drive, Suite 2700 \\ Los Angeles, CA 90045 \\
U.S.A. \\ acrans@lmu.edu}

\bigskip

\noindent\textbf{Thomas M.~Fiore} received a B.S. in Mathematics and
a B.Phil. in German at the University of Pittsburgh. He completed
his Ph.D. in Mathematics at the University of Michigan in 2005 under
the direction of Igor Kriz. He is an N.S.F. Postdoctoral Fellow and
Dickson Instructor at the University of Chicago, and was a Profesor
Visitante at the Universitat Aut\`{o}noma de Barcelona during
2007-08. His research interests include algebraic topology, higher
category theory, and mathematical music theory.

\noindent\textit{Department of Mathematics \\ University of Chicago
\\ Chicago, IL 60637 \\ U.S.A. \\ and \\
Departament de Matem\`{a}tiques \\ Universitat Aut\`{o}noma de Barcelona \\
08193 Bellaterra (Barcelona)  \\ Spain \\ fiore@math.uchicago.edu}

\bigskip

\noindent\textbf{Ramon Satyendra} received his Ph.D. from the
University of Chicago in the History and Theory of Music. He is
currently an Associate Professor of Music Theory at the University
of Michigan and is on the editorial boards of {\it Int\'{e}gral} and
the {\it Journal of Mathematics and Music.}

\noindent\textit{School of Music, Theatre and Dance \\ University of
Michigan
\\ E.V. Moore Building \\ 1100 Baits Dr. \\
Ann Arbor, MI 48109-2085 \\ U.S.A. \\ ramsat@umich.edu}

\end{document}